\documentclass[11pt,reqno,a4paper]{amsart}
\usepackage[active]{srcltx}
\usepackage{amssymb,amsmath}
\usepackage{xcolor}
\usepackage[normalem]{ulem}
\def\erre{\mathbb{R}}
\def\m1{{I\!\!M}}
\newcommand{\rdue}{\erre^2}
\newcommand{\ci}{\mathbb{C}}

\newcommand{\lapl}{\Delta}

\renewcommand{\to}{\rightarrow}

\newcommand{\rife}[1]{(\ref{#1})}

\newcommand{\scp}{\scriptscriptstyle}
\newcommand{\dsp}{\displaystyle}

\newcommand{\finedim}{\hspace{\fill}$\square$}

\newcommand{\ii}{\infty}

\newcommand{\al}{\alpha}

\renewcommand{\th}{\vartheta}
\newcommand{\sg}{\sigma}

\newcommand{\om}{\Omega}
\newcommand{\lm}{\lambda}
\newcommand{\N}{\mathbb{N}}

\newcommand{\tgs}{\tilde{g}_s}

\newcommand{\tgz}{\tilde{g}_{\scp 0}}

\newcommand{\e}[1]{{\,\dsp e^{\dsp #1}}}

\newtheorem{theorem}{Theorem}[section]
\newtheorem{proposition}[theorem]{Proposition}
\newtheorem{lemma}[theorem]{Lemma}
\newtheorem{corollary}[theorem]{Corollary}
\newtheorem{remark}[theorem]{Remark}
\newtheorem{definition}[theorem]{Definition}
\newcommand{\brm}{\begin{remark}}
\newcommand{\erm}{\end{remark}}
\newcommand{\bdf}{\begin{definition}\rm}
\newcommand{\edf}{\end{definition}}
\newcommand{\bte}{\begin{theorem}}
\newcommand{\ete}{\end{theorem}}
\newcommand{\bpr}{\begin{proposition}}
\newcommand{\epr}{\end{proposition}}
\newcommand{\ble}{\begin{lemma}}
\newcommand{\ele}{\end{lemma}}
\newcommand{\bco}{\begin{corollary}}
\newcommand{\eco}{\end{corollary}}
\newcommand{\beq}{\begin{equation}}
\newcommand{\eeq}{\end{equation}}
\newcommand{\bdm}{\begin{displaymath}}
\newcommand{\edm}{\end{displaymath}}

\newcommand{\graf}[1]{\left\{\begin{array}{ll}#1\end{array}\right.}

\begin{document}
\numberwithin{equation}{section}
\parindent=0pt
\hfuzz=2pt
\frenchspacing

\title[Sharp Harnack-type inequalities]{Sharp Harnack-type inequalities for conformal metrics on Riemann surfaces with positive Gaussian curvature and
conical singularities}

\author[D. Bartolucci]{Daniele Bartolucci}

\author[P. Cosentino]{Paolo Cosentino}

\address{Daniele Bartolucci, Department of Mathematics, University of Rome {\it ``Tor Vergata"} \\  Via della ricerca scientifica n.1, 00133 Roma, Italy. }
\email{bartoluc@mat.uniroma2.it}

\address{Paolo Cosentino, Department of Mathematics, University of Rome {\it ``Tor Vergata"} \\  Via della ricerca scientifica n.1, 00133 Roma, Italy. }
\email{cosentino.research@gmail.com , cosentino@mat.uniroma2.it}


\thanks{ D.B. and P.C.
are partially supported by the MIUR Excellence Department Project MatMod@TOV
awarded to the Department of Mathematics, University of Rome ``Tor Vergata'' and by PRIN project 2022, ERC PE1\_11,
``{\em Variational and Analytical aspects of Geometric PDEs}'', by the E.P.G.P. Project sponsored by the University of Rome ``Tor Vergata'' and by INdAM-GNAMPA Project, CUP E53C25002010001.}

\begin{abstract} We solve a long standing open problem (G. Tarantello, PNLDE 72 (2008)) about 
Harnack-type inequalities for metrics on
Riemann surfaces with positive Gaussian curvature and conical singularities of positive order. We develop two distinct approaches. The first addresses the particular case of constant Gaussian curvature and relies on a refined representation formula for solutions of the singular Liouville equation, combined with a priori estimates for a class of Fuchsian equations. In the general case of strictly positive and Lipschitz continuous Gaussian curvature, we establish a refined estimate that enables us to recover oscillation control for solutions to the Gauss curvature equation, even in the most unfavorable scenario.   
\end{abstract}

\maketitle {\bf Keywords}: conformal metrics on Riemann surfaces, conical singularities, sup + inf
and sup$\times$inf inequalities, singular Liouville equation, Gauss equation.

\section{Introduction}
\setcounter{equation}{0}
A conformal metric $\tgz$ of class $C^2$ on a Riemann surface $S$ is defined in any
local coordinate system $z\in \om \subset \ci$ as follows,
$$
g_0(z)=\e{u(z)}|dz|^2,\quad z\in \om,
$$
where $u\in C^2(\om)$. Here and in the rest of this paper
$\om\subset \ci$ will always denote an open and bounded set.
The quantity
$$
\rho_0(z)=\e{u(z)},
$$
is said to be the (local) conformal factor. Let $P_0\in S$, a conformal metric $\tgs$ on $S\setminus\{P_0\}$ is said to have a conical singularity
of order $\al\in(-1,0)\cup (0,+\ii)$
(or of angle $\th_\al=2\pi(1+\al)$) at  $P_0\in S$ if there exist $\om\subset \ci$ such that $0\in\om$ and 
local coordinates $z(P)\in \om\subset \ci$  and
$u\in C^{0}(\om)\cap C^{2}(\om\setminus \{0\})$ such that $z(P_0)=0$ and
$$
g_s(z)=|z|^{2\al} \e{u}|dz|^2,\quad z\in\om,
$$
where $g_s$ is the local expression of $\tgs$. We will denote by
$$
\rho_s(z)=|z|^{2\al} \e{u(z)}
$$
the singular (local) conformal factor and by
$$
\rho_{s,0}=|z|^{-2\al}\rho_s(z)=\e{u(z)},
$$
its regular part.
It is well known that the Gaussian curvature $K$ induced by $\tgs$ is defined  in local coordinates by the Gauss equation
\beq\label{Li}
-\Delta u = 2|z|^{2\al}K\e{u}\;\;\;\;\mbox{in}\;\;\;\;\;\om.
\eeq
Assuming $K$ to be at least locally Lipschitz, if the metric has no singularities, then $u$ satisfies \rife{Li} with $\al=0$ and is therefore of class $C^{2,\sg}_{\rm loc}(\om)$, by standard elliptic theory. In particular the conformal factor of
a metric on a Riemann surface of constant Gaussian curvature is always strictly positive and real analytic.
This well known geometrical interpretation was the motivation for the analysis of \rife{Li} initiated by Liouville
\cite{Lio} and Picard \cite{Pic}. In particular in \cite{Lio} a well known representation formula was derived for solutions of
\rife{Li} with $\al=0$ and $K\equiv \kappa\in \mathbb{R}$. A generalization of the Liouville formula to the case $\al\neq 0$ has been derived in \cite{CW}, \cite{bt2}.
We also refer to \cite{T2} for an account of the relevance of equation \rife{Li} in gauge field theories.\\
A systematic study of the regularity theory and
concentration compactness results for Liouville-type equations has been initiated in
a pioneering work by Brezis and Merle, see \cite{bm}.
On the basis of a conjecture in \cite{bm} and by using the Liouville formula, it was shown (among other things) in \cite{S}  that
the following holds:\\
Fix $\al=0$ and $K\equiv \kappa>0$ in $\om$. For any relatively compact subset $\om_0\subset\subset \om$, there exists a universal constant
$C_\om>0$ depending only on the distance of $\om_0$ from $\partial \om$, such that
\beq\label{1.i}
\sup\limits_{\om_0}{u}+\inf\limits_{\om} u\leq C_\om-2\log{(2\kappa)},
\eeq
for any distributional solution of \rife{Li}. Concerning the ``regular" case, that is $\al=0$, and if $K$ satisfies
\beq\label{3}
0<a\leq K\leq b<+\ii,
\eeq
and 
\beq\label{3.1}
\|\nabla K\|_{L^{\ii}(\om)}\leq L,
\eeq
the inequality \eqref{1.i} was established in \cite{bls}, with $C_\om-2\log(2\kappa)$ replaced by some constant $C_{\om,a,b,L}$, depending also on $a,b,L$.
A $``\sup+ C_1 \inf"$ inequality for
$K$ satisfying \rife{3} and just assumed to be bounded, where $C_\om$ and $C_1$ depend also on $a$ and $b$,
was established in \cite{S}. To the best of our knowledge, the most general result at hand so far has been obtained in \cite{cl},
where the optimal $C_1=\sqrt{\frac{b}{a}}$ was derived in terms of $a$
and $b$ for $K$ just assumed to satisfy \eqref{3}, as well as that $C_1=1$ for any $K$ sharing a uniform logarithmic modulus of continuity. It is worth recalling that a major breakthrough about estimates for Liouville-type equations was obtained in \cite{yy}, 
showing that, assuming a uniform bound on the boundary oscillation 
of solutions of \eqref{Li}, then, in the standard framework of the blow up phenomenon,  one can prove that $\sup\limits_{\om_0}{u}+\inf\limits_{\om} u\geq C$, 
provided that $\om_0$ contains an interior blow up point, yielding refined pointwise estimates for blow up solutions. This was in turn the starting point of several refined results, see for example \cite{bart-5, bart-6,cl-3,dwz,Lin-0,T4,T5,wz1,zhang1,zhang2} and references quoted therein.\\

As a Corollary of the results in  \cite{S} and \cite{bls} we have
\begin{corollary}[\cite{S}, \cite{bls}]
Let $\rho_0:\om\mapsto (0,+\ii)$ be the local conformal factor of a metric of  positive Gaussian curvature $K$ satisfying \eqref{3}, \eqref{3.1}.\\
For any relatively compact subset $\om_0\Subset \om$, there exists a positive constant
$C>1$ depending on the distance of $\om_0$ from $\partial \om$ and $a,b,L$, such that
\beq\label{2}
\left(\sup\limits_{\om_0}\rho_0\right)\cdot \left(\inf\limits_{\om}\rho_0\right)\leq C.
\eeq
\end{corollary}

\bigskip

This is the two dimensional analogue of the $``\sup \, \times \, \inf"$ inequalities for solutions
of the Yamabe equation in higher dimensions see \cite{Shoen} and \cite{cl-1}, \cite{Li-1}, \cite{Li-2}, \cite{LiZ0}, \cite{LiZ}.\\

The ``singular" case, that is $\al\neq 0$, is more involved. We always assume in this situation that $0\in\om$.
The generalization of the results in \cite{S},\cite{bls}, \cite{cl} is known in case
$\al\in(-1,0)$, see \cite{B1}, \cite{B2} and more recently \cite{Cos}.\\ 
At least to our knowledge, the case $\al>0$, first raised in \cite{bcct}, was still open to date and posed as a challenging open problem in \cite{T2} (see \cite{T2} \S 5.5.3). The issue is not trivial because the weight function $|z|^{2\al}K$ is not bounded from below away from zero, whence it fails to satisfy the left hand side inequality in \eqref{3}. Interestingly enough, it has been recently shown that, in some closely related situations, such as those involving collapsing singularities or regularization of singular sources,
a uniform upper bound of $\inf_{\om}u$ in terms of $-\sup_{\om_0}u$ does not hold in general (\cite{llty}). This is related to the phenomenon of blow up without mass concentration, see \cite{llty} and references quoted therein. However, the ``sup+inf" inequality was obtained in \cite{bcct}, under additional assumptions (in the same spirit of \cite{yy}) on the boundary oscillation of $u$. The best result at hand so far, which does not make use of other
assumptions, seems to be the sharp $u(0)+\inf\limits_{\om} u\leq C$ inequality obtained in \cite{T1}, see also \cite{T3}, \cite{T2}. Other partial results on the case $\al>0$ have been obtained in \cite{B3}, \cite{B4} and more recently \cite{BCW0}. We solve this problem here, providing two distinct proofs, one of which is about the case of independent interest where the curvature is constant, $K\equiv \kappa>0$ in $\om$.\\

Our first result is a sharp generalization of the estimate in \cite{S}. Indeed, we have the following:
\bte\label{supinf-sing}
Let $\al\in(0,+\ii)$ and $K\equiv \kappa>0$ in $\om$.
For any relatively compact subset $\om_0\Subset \om$, there exists a universal constant
$C_{\om,\al}>0$ depending only on $\al$ and on the distance of $\om_0$ from $\partial \om$ such that
\beq\label{1}
\sup\limits_{\om_0}{u}+\inf\limits_{\om} u\leq C_{\om,\al}-2\log(2\kappa),
\eeq
for any distributional solution of {\rm \rife{Li}}.
\ete
It is well known that the coefficient of the infimum in \eqref{1} is sharp, as can be seen from the classical explicit solutions, 
\[
 u_n(z)=\log\left(\frac{4(1+\al)^2n^2}{(1+n^2|z|^{2(1+\al)}\bigr)^2}\right),
\]
which solve \eqref{Li} in $\{|z|<1\}$ with $K\equiv \kappa=1$ and, as $n\to\infty$, $\sup_{|z|<{1/2}}u_n=2\log n+O(1)$ while $\inf_{|z|<1}u_n=-2\log n+O(1)$.

\bigskip

\bigskip

Here and in the rest of this paper $D_r$ will denote the disk of radius $r>0$, $D_r=\{z\in\ci\,|\,|z|<r\}$,
and $D^{*}_r=\{z\in\ci\,|\,0<|z|<r\}$, the corresponding punctured disk. We also set $D=D_1$.\\

The proof of Theorem \ref{supinf-sing} requires an improvement of the method in \cite{B4}, via a generalization of an argument in \cite{St} which,
as observed in \cite{GKR} (see also \cite{FKR}), yields a refinement of the estimate in \cite{S} about the regular case. To highlight the analogy with the singular setting, we briefly recall how the argument applies when $\alpha=0$. Let $w$ be a solution of \eqref{Li} with $\alpha=0$ and $K\equiv\kappa=1/2$ in the unit disk $D$. By the classical Liouville representation formula (\cite{Lio}), there exists a locally univalent meromorphic function $g$ in $D$ such that 
\[
w(z)=\log\left(\frac{8|g'|^2}{(1+|g|^2)^2}\right).
\]
Let us define the spherical derivative of $g$,
\[
g^{\#}(z):=\frac{|g'(z)|}{1+|g(z)|^2}.
\]
Therefore, $e^{\displaystyle w(z)}=8(g^\#(z))^2$
which implies that
\begin{equation}
    \label{estimatebelow}
g^\#\geq \frac{e^{\displaystyle \tfrac m2}}{\sqrt 8}:=\epsilon,
\end{equation}
with $m:=\inf_{D}w$. We may assume $m>-\infty$, since otherwise the desired inequality is trivial. Following the argument in \cite{St}, the locally univalent meromorphic function $g$ can be represented as $g=w_1/w_2$, where $w_1$ and $w_2$ are holomorphic in $D$ and satisfy the normalized Wronskian identity, $w_1'w_2-w_1w_2'=1$. Consequently, it is easy to see that
\begin{equation}\label{ff}
g^\#=\frac{1}{|w_1|^2+|w_2|^2}.
\end{equation}
By \eqref{estimatebelow}, $\|w_1\|_\infty,\|w_2\|_\infty\leq \tfrac{1}{\sqrt{\epsilon}}$. On the other hand, for $|z|\leq 1/2$, by using the Cauchy formula, we estimate the first derivative of $w_i$, namely $|w_i'(z)|\leq\tfrac{2}{\sqrt{\epsilon}}$ and by the Wronskian identity, we obtain
\[
1\leq \tfrac {2}{\sqrt{\epsilon}}(|w_1|+|w_2|)\leq \tfrac {2\sqrt 2}{\sqrt{\epsilon}}\sqrt{|w_1|^2+|w_2|^2}.
\]
Therefore, 
\[
\frac{e^{\displaystyle \tfrac{w(z)}{2}}}{\sqrt 8}=g^\#(z)\leq \frac{8}{\epsilon}=8\sqrt 8e^{\displaystyle-\tfrac{m}{2}},
\]
for every $|z|\leq 1/2$, which implies
\[
\sup_{D_{1/2}}w+\inf_D w\leq 2\log (64).
\]
Let us now consider the singular case $\alpha>0$. We first need a refinement of a Liouville-type representation formula obtained in \cite{CW} (see also \cite{bt2}) for solutions of \eqref{Li}.
\bte\label{thm:Liouville}
Any distributional solution $v$ of problem {\rm \rife{Li}} with $\al\in(0,+\ii)$, $K\equiv \kappa=1/2$ and $\om=D_4$
takes the following form,
\begin{equation}\label{eq:Lio1} v(z)=\log\left(\frac{8}{\bigl(|h(z)|^2+|z|^{2(1+\al)}|g(z)|^2\bigr)^2}\right),\quad z\in D_4,
\end{equation}
with $h$ and $g$ holomorphic in $D_4$ and satisfying,
\begin{equation}\label{eq:Lio2}
L_{\al}(z):=(1+\al)h(z)g(z)+z(h(z)g'(z)-h'(z)g(z))\equiv 1.
\end{equation}
In particular $(1+\al)h(0)g(0)=1$ and $h$, $g$ have no common zeroes.
\ete
We notice that, up to the factor $z^{\al}$, $L_\al$ is the Wronskian of the functions $z^{1+\al}g(z),h(z)$.
At this point, we apply the same argument as in the regular case described above, although unfortunately the Cauchy formula just implies $|z h'(z)|\leq 2|z|\|h\|_{L^{\ii}(D)}$ for any $|z|\leq 1/2$, which of course is  not enough to obtain
the needed generalization of the result in \cite{S}. The solution of this problem comes from the fact that,
taking the derivative of $L_\al$, gives that $h$ satisfies a Fuchsian equation with a peculiar structure, which allows a rather delicate improvement
of the Cauchy estimate, see Lemma \ref{lem:unifest} below.

\bigskip

Concerning the general case, we have the following sharp generalization  of the result in \cite{bls}, 
\begin{theorem}\label{thm:Kgen}
Let $\alpha\in(0,+\infty)$ and assume that $K$ satisfies \eqref{3}, \eqref{3.1}.
For any relatively compact subset $\om_0\subset\subset \om$, there exists
$C>0$, depending only on $\al$, $a$, $b$, $L$ and on the distance of $\om_0$ from $\partial \om$, such that
\beq\label{sup+infwithK}
\sup\limits_{\om_0}{u}+\inf\limits_{\om} u\leq C,
\eeq
for any distributional solution of 
\eqref{Li}.
\end{theorem}

The proof relies on the introduction of a new scaling distance, which allows one to recover a sort of control on the boundary oscillation of the solutions, as we shortly describe hereafter. By contradiction, using the $u(0)+\inf u$ 
inequality in \cite{T1} and the ``regular" $\sup+\inf$ inequality in \cite{bls}, we are reduced to analyzing the case where $\sg_n:=u_n(x_n)+m_n\to +\ii$, with $m_n=\inf\limits_{\om} u_n$, $u_n(x_n)=\sup\limits_{\om_0} u_n$, $|x_n|\to 0$,  
$|x_n|^{2(1+\al)}\exp{(u_n(x_n))}\to +\ii$ and 
$\sg_n\leq -4(1+\al)\log(|x_n|) +C$. 
Of course, these facts are well known and 
yield no contradiction, since $|x_n|\to 0$.\\
The novelty of the argument is 
the introduction of new scales,  $R_n^{2(1+\al)}=\exp{(\sg_n/8)}\to +\ii$ and $s_n=|x_n|R_n\to 0$, which allows 
one to solve some rather subtle indeterminate forms. First of all, by a nice trick based on the Green representation formula and $\sg_n\leq -4(1+\al)\log(|x_n|) +C$, it can be shown that the mass contained in $D_{s_n}$ is strictly less than $8\pi(1+\al)$. This is somehow inspired by ideas in \cite{cl}, see also \cite{Cos}. On the other side, a careful use of the Harnack inequality and elliptic estimates shows that the oscillation of the term in the Green representation formula due to the exterior of $D_{s_n}$ is small. At this point, the crux of the proof is to show that, scaling at distance $|x_n|$ and using the information gathered so far at scale $s_n=|x_n|R_n>>|x_n|$, one has a uniform control of the oscillations of the solutions, far away from the possibly non empty blow up set. This is done by careful elliptic estimates, see Lemma \ref{lemmaconvergenzariscalate}. 
Once this is done, we find ourselves exactly in the situation already discussed in \cite{bcct} via a moving plane argument.

\bigskip

Concerning our original motivation, namely the ``$\sup\times\inf"$ inequalities, there is
no direct nontrivial analogue of \rife{2} for $\rho_s$. Indeed, as far as $\alpha>0$, then $\inf_{\Omega}\rho_s=0$ if  $0\in\Omega_0$, 
whence the corresponding
product vanishes.
A meaningful generalization of \rife{2} to the singular case is obtained by looking at $\rho_{s,0}$,
the regular part of the local conformal factor.\\ 
As an immediate Corollary of our results we have,
\begin{corollary}
Let $\rho_{s,0}:\om\mapsto (0,+\ii)$ be the regular part of the local conformal factor of a metric of
Gaussian curvature $K$ satisfying \eqref{3}, \eqref{3.1}, on a Riemann surface $S$ with a conical singularity of order
$\al\in(0,+\ii)$.\\
For any relatively compact subset $\om_0\Subset \om$ there exists a positive constant
$C>1$ depending only on $\al,a,b,L$ and on the distance of $\om_0$ from $\partial \om$, such that
\[
\left(\sup\limits_{\om_0}\rho_{s,0}\right)\cdot \left(\inf\limits_{\om}\rho_{s,0}\right)\leq C.
\]
\end{corollary}

This paper is organized as follows. In section \ref{sec2}, we give the proof of the representation formula. In section \ref{sec3} we derive a priori estimates for a Fuchsian equation, while in section \ref{sec4} we give the proof of Theorem \ref{supinf-sing}. In section \ref{sec5} we prove Theorem \ref{thm:Kgen}.

\bigskip
\bigskip

\section{The proof of Theorem \ref{thm:Liouville}}\label{sec2}

\emph{Proof of Theorem \ref{thm:Liouville}.}\\
Put $\beta=1+\alpha$ and let $v$ be any distributional
solution of \eqref{Li} with $K\equiv 1/2$. We prove \eqref{eq:Lio1} and \eqref{eq:Lio2}
for $v$ on the full disk $D_4$. Since $v$ is a distributional
solution, it is understood that $|z|^{2\alpha}e^v\in L^1_{\rm loc}(D_4)$, whence it is well
known (\cite{bm}) that $v\in C^{2}(D_4)$.\\
By Theorem~3 in \cite{CW} for $u=v+2\alpha\log|z|$ on $D_4^*$, there exists a meromorphic and locally univalent
developing map $F$ on the universal cover of $D_4^*$, which is $\widetilde D_4^*=\{(r,\theta)\in (0,4)\times \mathbb{R}\}$. In the local coordinate $z=re^{i\theta}$,
\[
 F^\#(z)=\frac{|F'(z)|}{1+|F(z)|^2}
     =\frac{1}{\sqrt{8}}|z|^{\beta-1}e^{v/2},
\]
where the spherical derivative at a pole is computed using $1/F$.
In particular, 
$$|F^\#(z)|\le C |z|^{\beta-1}$$
near zero. Let $T(r,\theta)=(r,\theta+2\pi)$ be the deck transformation corresponding to a positively oriented circuit around the puncture. Since $F$ and $F\circ T$ induce the same spherical metric, they differ by a local isometry $M$ of the sphere. Consequently, $F\circ T=M\circ F$. In fact, it is well known that a local isometry of the sphere is the restriction of a
unique global isometry (see for example Corollary 6.4, p. 256 in \cite{KN}), which implies that  $M\in\mbox{O}(3)$. On the other side, since the action is obtained locally by composition
of holomorphic maps, it is holomorphic itself, whence preserves the orientation. Therefore $M\in S\mbox{O}(3)\simeq PSU(2)$,
which takes the general form 
$$M(w)=\frac{aw+b}{-\overline{b}w+\overline{a}},\qquad |a|^2+|b|^2=1.$$
Let $\gamma(t)=t e^{i\theta_0}$, $0<s\leq t\leq r<4$ be any radial segment and
$\widetilde{\gamma}$ its lift to the universal cover, then the length of $F \circ \tilde{\gamma}$ with respect to the metric $ds_{sph}=|dw|/(1+|w|^2)$
is at most 
\begin{equation}
    \label{radialestimate}
    L_{sph}(F\circ\tilde\gamma)\leq C\frac{(r^\beta-s^\beta)}{\beta}.
\end{equation}
Thus, for each fixed $\theta$, $F(\tilde{\gamma})$ has a spherical limit along each
lifted ray as $r\to 0$, say $p_\theta$. On the other side, any lifted arc $\gamma(t)=r e^{it}$, $\theta_1 \leq t \leq \theta_2$,
between two fixed angles, has image of length at most
$Cr^\beta|\theta_2-\theta_1|$. Therefore all the spherical limits $p_\theta$ as $r\to 0$ are the same point. We denote their common value by $p$.
It is not difficult to see that
in particular $M(p)=p$. Indeed, denoting the general coordinate on the covering $(r,\theta), r\in(0,4), \theta\in\mathbb{R}$,
after one turn we have
\[
F(r,\theta+2\pi)=M(F(r,\theta)).
\]
Thus, taking the limit as $r \to 0$ on the left gives $p$, while by the continuity of $M$, the limit on the right gives $M(p)$.
At this point we compose with a spherical isometry sending $p$ to $0$. Therefore, denoting once again the new isometry with $M$, we have $M(0)=0$, 
which implies $b=0$ and $|a|=1$, that is, the monodromy map
is $w\mapsto e^{2\pi i\nu}w$, for some $0\le\nu<1$.
By using \eqref{radialestimate} and taking the limit as $s\to0$, we deduce that
\[
d_{sph}(F(r,\theta),0)\leq\tfrac{C}{\beta}r^{\beta},
\]
uniformly in $\theta$. Since the metric is $|dw|/(1+|w|^2)$, the spherical distance is $d_{sph}(w,0)=\arctan|w|$, therefore $$\arctan|F(r,\theta)|\leq\tfrac{C}{\beta}r^{\beta},
\,\,\text{uniformly in}\,\, \theta,$$ 
which implies that $|F|=O(r^\beta)$, uniformly in $ \theta$.\\
It is easy to verify that $\Phi (z):=z^{-\nu}F(z)$ descends to a meromorphic function on
$D_4^*$, which is holomorphic near zero and $\Phi(z)=O(r^{\beta-\nu})$.
Since $\beta>1$, it extends with a zero of order $k\ge1$. Therefore,
$F=z^{k+\nu}\psi$, with $\psi$ a single-valued and meromorphic function
on $D_4$, holomorphic at zero and $\psi(0)\ne0$. Comparing
\[
 F^\#=(k+\nu)|\psi(0)|\,r^{k+\nu-1}(1+o(1))
 \quad\hbox{and}\quad
 F^\#=8^{-1/2}e^{v(0)/2}r^{\beta-1}(1+o(1))
\]
gives $$k+\nu=\beta.$$ 

Set $B=\beta\psi+z\psi'$. Then $F'=z^{\beta-1}B$ on every local
branch. Local univalence gives $B\ne0$ at every regular point of
$\psi$ in $D_4^*$ and $B(0)=\beta\psi(0)\ne0$. Let $a$ be a pole of $\psi$, we have that $q(z)=1/F(z)$ is a local coordinate 
there, hence in particular $q(a)=0$ and $q'(a)\neq  0$. Therefore we can write $1/F(z)=q'(a)(z-a)+O((z-a)^2)$, that is $F(z)$ and consequently $\psi(z)$ 
have a simple pole. Therefore any pole $a$ of $\psi$ is simple and writing $\psi=c/(z-a)+O(1)$, $c\neq 0$, gives
\[
 B=-\frac{ac}{(z-a)^2}+O\bigl((z-a)^{-1}\bigr),\qquad ac\ne0 .
\]
Consequently $1/B$ is holomorphic on $D_4$ and all its zeros are
double, since,
$$
\frac{1}{B(z)}=-\frac{(z-a)^2}{ac}\big(1+O\bigl((z-a)\bigr)\big),\qquad ac\ne0 .
$$ \\
On a disk, a holomorphic function with zeros
of even order has a holomorphic square root. Choose such a root $h$
of $1/B$ and put $g=h\psi$. At each pole of $\psi$ the simple zero
of $h$ cancels the pole, so $g$ is holomorphic on $D_4$.

Away from the poles, since $h^2=1/B$, we have
\[
 \beta hg+z(hg'-h'g)=h^2(\beta\psi+z\psi')=1 .
\]
Clearly the identity extends over $D_4$. In particular, it implies that $h$ and $g$ have no
common zero and that $\beta h(0)g(0)=1$. Moreover,
\[
 |h|^2+|z|^{2\beta}|g|^2
 =\frac{1+|F|^2}{|B|}.
\]
Substitution into the developing-map formula gives
\[
 e^{\displaystyle v}=\frac{8}{\bigl(|h|^2+|z|^{2\beta}|g|^2\bigr)^2}.
\]
Both sides are continuous on $D_4$, including the origin and the
poles of $\psi$. This proves the two identities.
\finedim
\bigskip
\bigskip

\section{A priori estimates for a Fuchsian equation} \label{sec3}
We will need the following a priori estimate.
\begin{lemma}\label{lem:unifest} Fix $\beta>1$, there exist $C_\beta>0$ such that if $h$ is holomorphic and bounded in $D$
and for some $q$ holomorphic in $D$ it holds,
\begin{equation*}
 zh''(z)-(\beta-1)h'(z)=q(z)h(z),
\end{equation*}
then
\begin{equation}\label{eq:lemma-bound}
 {|zh'(z)|\leq C_\beta\left(|h(z)|+|z|^\beta
                      \|{h}\|_{L^\infty(D)}\right),
             \qquad |z|\leq\tfrac12.}
\end{equation}
\end{lemma}
\proof
Since the estimate is homogeneous in $h$, and the case $h\equiv0$ is trivial, we may assume without loss of generality that
$\|h\|_{L^\infty(D)}=1$.
Argue by contradiction and assume that there exist
$h_n,q_n$ and $\eta_n\in D_{1/2}\setminus\{0\}$ such that $\|h_n\|_{L^\infty(D)}=1$ and
\begin{equation}\label{eq:failure}
 \frac{|\eta_n h_n'(\eta_n)|}
 {|h_n(\eta_n)|+|\eta_n|^\beta}\to+\ii.
\end{equation}
The Cauchy formula implies $|h_n'(z)|\leq2$ for any  $z\in D_{1/2}$, whence
\eqref{eq:failure} implies $\eta_n\to 0$. 
Let us set
\begin{equation}\label{eq:k}
 k=\lceil\beta\rceil-1,\qquad k<\beta\leq k+1,
\end{equation}
and write
\[
 h_n(z)=\sum_{j=0}^{\infty}a_{j,n}z^j,
 \qquad |a_{j,n}|\leq1,
\]
where the bound on the coefficients follows from the Cauchy estimates and the normalization
$\|h_n\|_{L^\infty(D)}=1$.
Clearly we have,
\begin{equation}\label{eq:derivative-tail}
 |\eta_n h_n'(\eta_n)|\leq C_k \max_{0\leq j\leq k}|a_{j,n}\eta_n^j|+C_k |\eta_n|^{k+1},
\end{equation}
for some $C_k\geq 1$. It follows from \eqref{eq:failure} that
$|\eta_n h_n'(\eta_n)|/|\eta_n|^\beta\to+\infty$ and, since $|\eta_n|^{k+1}\leq |\eta_n|^\beta$,
we see from \eqref{eq:derivative-tail} that
\begin{equation}\label{eq:sn-large}
 \frac{\max_{0\leq j\leq k}|a_{j,n}\eta_n^j|}{|\eta_n|^\beta}\to+\ii.
\end{equation}
In particular
\[
s_n:=\max_{0\leq j\leq k}|a_{j,n}\eta_n^j|>0
\]
for $n$ large and we can set
\begin{equation}\label{eq:rescaled}
 F_n(w)=\frac{h_n(\eta_n w)}{s_n}=\sum_{j=0}^{+\ii}b_{j,n}w^j,
 \qquad |w|<1/|\eta_n|,
\end{equation}
where $|b_{j,n}|\leq 1$ for any $j\leq k$ and at least one among these coefficients satisfies $|b_{j,n}|= 1$.
For each $R>0$ fixed, the remainder satisfies,
\[
 \sup_{|w|\leq R}\left|F_n(w)-\sum_{j=0}^k
             \frac{a_{j,n}\eta_n^j}{s_n}w^j\right|
 \leq C_{R,k}\frac{|\eta_n|^{k+1}}{s_n}\longrightarrow0.
\]
Therefore, passing to a subsequence, we have $b_{j,n}\to b_j$, with
$\max_{0\leq j\leq k}|b_j|=1$. Hence,
\begin{equation}\label{eq:P-limit}
 F_n\to P\quad\text{locally uniformly in }\mathbb{C},
 \qquad P\not\equiv0,\quad \deg P\leq k<\beta,
\end{equation}
with local uniform convergence for the derivatives of $F_n$ as well.

Moreover \eqref{eq:derivative-tail} and \eqref{eq:sn-large} give
$$|\eta_n h_n'(\eta_n)|/s_n\leq C_k+o(1),$$ while
from \eqref{eq:failure} we have,
$$|h_n(\eta_n)|/|\eta_n h_n'(\eta_n)|\to0.$$
Therefore,
\begin{equation}\label{eq:P-one}
 P(1)=\lim_n h_n(\eta_n)/s_n=0.
\end{equation}
On the other side, after rescaling, we deduce from our differential equation that,
\[
 wF_n''-(\beta-1)F_n'=H_nF_n,
 \qquad H_n(w)=\eta_n q_n(\eta_n w).
\]
Remark that each $H_n$ is holomorphic on  $D_{1/|\eta_n|}$.
For any $R$ large enough, $P$ has no zeroes on $|w|=R$, hence
\[
 H_n=\frac{wF_n''-(\beta-1)F_n'}{F_n}
 \longrightarrow\frac{wP''-(\beta-1)P'}P
\]
uniformly on $\partial D_R$. By the maximum principle we have $|H_n|\leq C_R$ on $\overline D_R$, implying that
the limiting quotient may admit only removable singularities in $D_R$.
Therefore, choosing larger and larger disks, we deduce that
\begin{equation}\label{eq:entire-rational}
 \frac{wP''-(\beta-1)P'}P
\end{equation}
is a rational entire function, whose numerator has degree strictly less than the degree of $P$. However, because of  \eqref{eq:P-one}, $\deg P\geq 1$
whence \eqref{eq:entire-rational} must vanish identically,
\begin{equation}\label{eq:Euler-P}
 wP''-(\beta-1)P'=0.
\end{equation}
At this point, writing $P(w)=\sum_{j=0}^k p_jw^j$, we deduce that,
\[
 \sum_{j=1}^k j(j-\beta)p_jw^{j-1}=0,
\]
but since $1\leq j\leq k<\beta$, all the coefficients $p_j$, $j\geq1$, must vanish as well.
Therefore $P$ should be a non vanishing constant, which again contradicts \eqref{eq:P-one}. This contradiction proves \eqref{eq:lemma-bound}.\\
Remark that if $\beta\in\N$, the Euler equation $wP''-(\beta-1)P'=0$ admits the solution $P=w^\beta$. However by definition in \eqref{eq:k} we have
$k<\beta$ which rules out that solution.
\finedim

\bigskip
\bigskip

\section{The proof of Theorem \ref{supinf-sing}}\label{sec4}
In this section we prove Theorem \ref{supinf-sing}.\\
\emph{The proof of Theorem \ref{supinf-sing}.}\\
Let $u$ be a solution of \eqref{Li} in $\om$, the function $u(rx)+2(1+\al)\log(r)$ solves the same equation in $\frac{1}{r}\om$, whence there is no loss of
generality in assuming $D_4\Subset \om$. In view of the regular ``sup+inf" inequality (\cite{S},\cite{bls}), it is enough to prove that any $v=u+\log(2\kappa)$ satisfies
\begin{equation}\label{est0}
\sup_{D_{\frac12}}{v}+\inf\limits_{\om}v\leq \sup_{D_{\frac12}}{v}+\inf_{D_4}v\leq C,
\end{equation}
for some $C$ that depends only on $\al$. Remark that $v$ is a solution of \eqref{Li} with $K\equiv 1/2$. Let
\[
m=\inf_{D_4}v,\quad \delta^4_m=e^{-m},
\]
and assume without loss of generality that $m>-\ii$. Let $h$ and $g$ be defined as in Theorem \ref{thm:Liouville}, then by the maximum 
principle we have,
\begin{equation}\label{est1}
\|h\|_{L^{\ii}(D)}\leq \sqrt[4]{8}\delta_m,\quad  \|g\|_{L^{\ii}(D)}\leq \sqrt[4]{8}\delta_m, \quad \|g'\|_{L^{\ii}(D_{\frac12})}\leq 2\sqrt[4]{8}\delta_m
\end{equation}
where we used the Cauchy formula for $g'$. Let $\beta=1+\al$, taking the derivative of \eqref{eq:Lio2} we have
\[
0=(\beta-1)h'g+(\beta+1)hg'+z(hg''-h''g),
\]
implying that
\[
g[zh''-(\beta-1)h']=h[zg''+(\beta+1)g'].
\]
Since $h$ and $g$ have no common zeroes, we see that
$$
q(z)=\frac{zh''(z)-(\beta-1)h'(z)}{h(z)},
$$
is holomorphic in $D_2$. Indeed, it is clearly holomorphic as far as $h(z)\neq 0$, while near a zero of $h$ we have
$$
q(z)=\frac{zg''(z)+(\beta+1)g'(z)}{g(z)},
$$
which is holomorphic since $g$ cannot be zero and the two expressions coincide away from the zeroes of $h$ or $g$. Therefore we deduce that
\[
zh''(z)-(\beta-1)h'(z)=q(z)h(z), \quad z\in D_2,
\]
and we can apply Lemma \ref{lem:unifest}, which gives, 
\[
|zh'(z)|\leq C_{\beta}(|h(z)|+|z|^\beta\|h\|_{L^{\ii}(D)})\leq C_{\beta}(|h(z)|+|z|^\beta\sqrt[4]{8}\delta_m),\quad |z|\leq \frac12,
\]
where we used \eqref{est1}. Recall that $C_{\beta}$ does not depend on $q$. We use this estimate together with \eqref{est1} in $1=L_\al(z)$, with $|z|\leq 1/2$, which gives,
\[
1\leq \beta |h||g|+|z||h||g'|+|z||h'||g|\leq \beta|h| \sqrt[4]{8}\delta_m + |h|\sqrt[4]{8}\delta_m + C_{\beta}(|h|+|z|^\beta\sqrt[4]{8}\delta_m)|g|\leq
\]
\[
(\beta +1 + C_\beta) |h| \sqrt[4]{8}\delta_m + C_{\beta}|z|^\beta\sqrt[4]{8}\delta_m |g|\leq \overline{C}_\beta (|h|+|z|^\beta|g|)\sqrt[4]{8}\delta_m,
\]
that is,
\[
1\leq 2\overline{C}_\beta^2 (|h|^2+|z|^{2(1+\al)}|g|^2)\sqrt{8} e^{-m/2}.
\]
Therefore,
\[
1\leq 16\overline{C}_\beta^2 e^{-(v(z)+m)/2},\qquad\text{for}\,\,|z|\leq\tfrac{1}{2}.
\]
which readily implies 
\[
v(z)+m\leq 2\log(16\bar C^2_\beta),\qquad\text{for}\,\,|z|\leq\tfrac{1}{2}
\]
and gives our thesis \eqref{est0}. \finedim

\bigskip

\section{The Proof of Theorem \ref{thm:Kgen}}\label{sec5}

\begin{proof}
We briefly show that it is enough to prove that
\begin{align}
    \label{sup+infreduced}
\sup_{\bar D_1}u+\inf\limits_{D_4} u\leq C.
\end{align}

Indeed, the function $u(rx)+2(1+\al)\log(r)$ solves an equation of the same form in $\frac{1}{r}\om$, whence there is no loss of
generality in assuming $D_4\Subset \om$.  Moreover, let $\delta={\rm dist}(\Omega_0,\mathbb{R}^2\setminus\Omega)>0$. For $x\in \Omega_0$, with $|x|\leq\delta/16$, the disk $D_{\delta/2}$
is contained in $\Omega$. Let
$R=\delta/8$ and
\[
 \widehat u(y)=u(Ry)+2(1+\alpha)\log R,
 \qquad \widehat K(y)=K(Ry),\qquad \|\nabla\widehat K\|_\infty\leq LR.
\]
Since $x/R\in D_1$, \eqref{sup+infreduced} implies that
\[
 u(x)+\inf_\Omega u\leq C(\alpha,a,b,LR)-4(1+\alpha)\log R.
\]
For $|x|>\delta/16$, let
$\rho=\min\{\delta/4,|x|/4\}\geq\delta/64$.
Then, $B_{2\rho}(x)$ is contained in $\Omega$ and $0\notin B_{2\rho}(x)$.
Therefore,
\[
 w(y)=u(x+\rho y)+2\log\rho+2\alpha\log|x|,\qquad |y|<2,
\]
satisfies a regular Liouville type equation with potential
\[
 \left|\frac{x}{|x|}+\frac\rho{|x|}y\right|^{2\alpha}K(x+\rho y),
\]
which satisfies the hypotheses of \cite{bls}, thus
\[
 u(x)+\inf_\Omega u
 \leq C-4\log\rho-4\alpha\log|x|
 \leq C(\alpha,a,b,L,\delta).
\]
Taking the supremum in $\Omega_0$, we obtain the desired inequality \eqref{sup+infwithK}.

\bigskip

Therefore, we are left to show the validity of \eqref{sup+infreduced}. We argue by contradiction and suppose that there exist sequences $\tilde u_n$ and $K_n$ satisfying
\[
-\lapl\tilde u_n=2|z|^{2\alpha}K_ne^{\displaystyle \tilde u_n}\quad\text{ in }\quad D_4,
\] 
\[
0<a\leq K_n\leq b<+\infty \qquad\text{and}\qquad \|\nabla K_n\|_{L^{\infty}(D_4)}\leq L,
\]
such that
\[
\sup_{\bar D_1}\tilde u_n+\inf_{D_4}\tilde u_n\to+\infty.
\]
Let us set $u_n=\tilde u_n+\log2$. Therefore,
\[
-\lapl u_n=|z|^{2\alpha}K_ne^{\displaystyle u_n}\quad\text{ in }\quad D_4
\]
and
\[
\sup_{\bar D_1} u_n+\inf_{D_4} u_n\to+\infty.
\]
Without loss of generality, we can assume that 
$$m_n:=\inf_{D_4}u_n>-\infty,\qquad \sup_{\bar D_1}u_n=u_n(x_n),$$
for some sequence $x_n$, $|x_n|\leq1$, and set
\begin{equation}
\label{contradiction}
\sigma_n:=u_n(x_n)+m_n\to+\infty.
\end{equation}
Since $m_n\leq u_n(x_n)$, $u_n(x_n)\to+\infty$ and, by applying the regular $\sup+\inf$ inequality (\cite{bls}), it is easy to see that $|x_n|\to0$. If, along a subsequence $|x_n|\equiv0$, then, by applying the ``$u(0)+\inf u"$ inequality (\cite{T1}), we obtain a contradiction. Moreover, setting $$\delta_n:=e^{\displaystyle -\tfrac{u_n(x_n)}{2(1+\alpha)}}\to0,$$
we show that 
\begin{equation}
    \label{casobuono}
    \frac{|x_n|}{\delta_n}\to+\infty.
\end{equation}
Otherwise, passing to a subsequence, we may assume $\tfrac{|x_n|}{\delta_n}\leq \hat C$. Let
\[
w_n(y)=u_n(\delta_n y)+2(1+\alpha)\log\delta_n, \quad  y\in B_{2\hat C}(0),
\]
which satisfies,
\[
\graf{-\lapl w_n=|y|^{2\alpha}K_n(\delta_ny)e^{\displaystyle w_n}:=f_n,\quad\text{ in }\quad B_{2\hat C}(0), \\
w_n\leq w_n(\tfrac{x_n}{\delta_n})=0.}
\]
Since $0\leq f_n\leq C$, by using a standard argument based on the Harnack inequality as applied to $-w_n\geq0$, we have $w_n(0)\geq -C$. In particular, this implies that $$u_n(0)\geq -2(1+\alpha)\log\delta_n-C=u_n(x_n)-C$$
and then, by using again the ``$u(0)+\inf u"$ inequality (\cite{T1}), we infer that
\[
\sigma_n\leq u_n(0)+m_n+C\leq \tilde C,
\]
which contradicts \eqref{contradiction}.\\

At this point, let 
\[
\tilde u_n(z)=u_n(x_n+|x_n|z)+2(1+\alpha)\log|x_n|,\qquad |z|\leq\tfrac12,
\]
which satisfies
\[
\graf{-\lapl \tilde u_n=\tilde V_n(z)e^{\displaystyle \tilde u_n},\qquad |z|\leq\tfrac12,\\
\tilde V_n(z)=|\tfrac{x_n}{|x_n|}+z|^{2\alpha} K_n(x_n+|x_n|z),}
\]
where $\tilde V_n$ satisfies in $D_{\frac12}$ the assumptions of the regular $\sup+\inf$ inequality (\cite{bls}). Therefore, we deduce that
\[
\tilde u_n(0)+\inf_{D_{1/2}}\tilde u_n\leq C_0
\]
which implies
\begin{equation}
    \label{boundfromabovesup+inf}
\sigma_n\leq u_n(x_n)+\inf_{B_{|x_n|/2}(x_n)}u_n\leq -4(1+\alpha)\log|x_n|+C_0.
\end{equation}
At this point, we define
\[
R_n:=\exp\left(\frac{\sigma_n}{16(1+\alpha)}\right)\quad\text{and}\quad s_n:=|x_n|R_n.
\]
Then, by \eqref{contradiction} and \eqref{boundfromabovesup+inf}, we have that 
$$R_n\to+\infty,\quad s_n\leq C|x_n|^{\frac34}\to0\quad \text{and}\quad \log(\tfrac2{s_n})\geq -\tfrac34\log|x_n|-C.$$

Now, let us consider the Green function of the Laplacian in $D_2$, that is, in complex notation,
\[
G(x,y)=\frac1{2\pi}\log\left(\frac{|4-\bar yx|}{2|x-y|}\right),\qquad x,y\in D_2
\]
and let us use the Green representation formula
\begin{equation}
    \label{Greenrepresentation}
u_n(x)=h_n(x)+\int_{D_2}G(x,y)|y|^{2\alpha}K_n(y)e^{\displaystyle u_n(y)}\,dy,
\end{equation}
with $h_n$ harmonic in $D_2$, $h_n=u_n$ on $\partial D_2$. By the maximum principle, $h_n\geq m_n$. We want to show that
\begin{equation}
    \label{massinDsn}
    \int_{D_{s_n}(0)}|x|^{2\alpha}K_ne^{\displaystyle u_n}\,dx<8\pi(1+\alpha).
\end{equation}
We argue by contradiction and suppose $\int_{D_{s_n}(0)}|x|^{2\alpha}K_ne^{\displaystyle u_n}\,dx\geq8\pi(1+\alpha).$ We observe that, if $|x-x_n|\leq \tfrac{|x_n|}{2}$, then $|x|\leq\tfrac32|x_n|=\frac32 \frac{s_n}{R_n}$, and then it is not difficult to see that for any $y\in D_{s_n}(0)$ we have,
\[
G(x,y)\geq\frac1{2\pi}\left(\log\frac1{s_n}-C\right),
\]
for some uniform constant $C$. Therefore, by \eqref{Greenrepresentation}, for every $x$ such that $|x-x_n|\leq\tfrac{|x_n|}{2}$, we have that
\begin{align*}
u_n(x)&=h_n(x)+\int_{D_2}G(x,y)|y|^{2\alpha}K_n(y)e^{\displaystyle u_n(y)}\,dy\\
 &\geq m_n+\int_{D_{s_n}(0)}G(x,y)|y|^{2\alpha}K_n(y)e^{\displaystyle u_n(y)}\,dy\\
 &\geq m_n+4(1+\alpha)\left(\log\frac1{s_n}-C\right).
\end{align*}
Therefore, by using \eqref{boundfromabovesup+inf}, we deduce that
\begin{align*}
\sigma_n&=u_n(x_n)+m_n\leq u_n(x_n)+\inf_{D_{\frac{|x_n|}{2}}(x_n)}u_n+4(1+\alpha)\log s_n+C\\
&\leq4(1+\alpha)\log\left(\frac{s_n}{|x_n|}\right)+C=4(1+\alpha)\log R_n+C\\
&=\frac14\sigma_n+C,
\end{align*}
which is a contradiction to \eqref{contradiction}. This fact proves \eqref{massinDsn}.\\
Let 
\[
H_n(x)=h_n(x)+P_n(x),\qquad P_n(x)=\int_{D_2\setminus D_{s_n}}G(x,y)|y|^{2\alpha}K_n(y)e^{\displaystyle u_n(y)}\,dy,
\]
then the formula \eqref{Greenrepresentation} becomes
\begin{equation}
    \label{Greenrepresentation2}
    u_n(x)=H_n(x)+\int_{D_{s_n}}G(x,y)|y|^{2\alpha}K_n(y)e^{\displaystyle u_n(y)}\,dy.
\end{equation}
We notice that $H_n$ is harmonic in $D_{s_n}(0)$ and, for $x\in D_{s_n}(0)$, $P_n\geq0$ and
\[
 u_n(x)\geq H_n(x)\geq m_n.
\]
Since $\tfrac{|x_n|}{s_n}\to0$, $D_{\frac{|x_n|}{2}}(x_n)\subset D_{s_n}(0)$ and, by applying the standard Harnack inequality to the function $H_n-m_n$, there exists a uniform constant $c>0$ such that
\[
c(H_n(x_n)-m_n)\leq \inf_{B_{\frac{|x_n|}{2}}(x_n)}(H_n-m_n).
\]
In particular, by using again \eqref{boundfromabovesup+inf}, we deduce that
\begin{align*}
\sigma_n+c(H_n(x_n)-m_n)&\leq \sigma_n +\inf_{B_{\frac{|x_n|}{2}}(x_n)}(H_n-m_n)\\
&=u_n(x_n)+\inf_{B_{\frac{|x_n|}{2}}(x_n)}H_n\\
&\leq u_n(x_n)+\inf_{B_{\frac{|x_n|}{2}}(x_n)}u_n\leq -4(1+\alpha)\log|x_n|+C_0,
\end{align*}
which, for $n$ large enough, implies that there exists a positive constant $C$ such that
\begin{equation}
    \label{controlHn}
    0\leq H_n(x_n)-m_n\leq -C\log|x_n|
\end{equation}
and in particular, since both $h_n(x_n)-m_n$ and $P_n$ are nonnegative, that
\begin{equation}
    \label{controlhnePn}
    0\leq h_n(x_n)-m_n\leq -C\log|x_n|,\quad\quad
    0\leq P_n(x_n)\leq -C\log|x_n|.
\end{equation}

At this point, we prove the following technical lemma

\begin{lemma}
    \label{lemmatecnico}
    There exists $C>0$ such that
    \begin{equation}
        \label{refinedgradientestimate}
    |\nabla_xG(x,y)|\leq\frac{C}{s\log(2/s)}G(x_0,y),
        \end{equation}
for $0<s<\tfrac14$, $|x|,|x_0|\leq s/2$ and $y\in D_2\setminus D_s$.
\end{lemma}
\begin{proof}
If $s\leq|y|\leq 1$, by differentiating the Green function, we immediately obtain
\[
|\nabla_xG(x,y)|\leq \frac{C}{|y|},\qquad G(x_0,y)\geq c\log\left(\frac{2}{|y|}\right).
\]
Moreover, we notice that the function $r\log(2/r)$ is concave in $[s,1]$ and for $s\leq\frac14$, has minimum in $r=s$. Therefore, 
\[
|y|\log(2/|y|)\geq s\log(2/s)
\]
and
\[
|\nabla_x G(x,y)|\leq\frac{C_1}{|y|\log(\frac{2}{|y|})}\log(\tfrac{2}{|y|})\leq \frac{C_2}{s\log(\frac{2}{s})}G(x_0,y),
\]
which is the thesis. If $1\leq|y|< 2$, the function $v_y(x)= G(x,y)$ is harmonic for $x\in D_{1/2}$ and, by standard Harnack and interior gradient estimates, we have that
\[
\sup_{D_{1/8}} |\nabla v_y|\leq C\,\sup_{D_{1/4}} v_y\leq C\,\inf_{D_{1/4}}v_y\leq C\, v_y(x_0)=C \,G(x_0,y).
\]
Since $s\log(2/s)<1$, we also have the desired result in this case.
\end{proof}

Let us fix $T>0$. By using standard interior gradient and Harnack estimates for the positive function $h_n-m_n$ and \eqref{controlhnePn}, we deduce that
\begin{align}
   \nonumber
    |x_n|\sup_{|x|\leq T|x_n|}|\nabla h_n|&\leq|x_n|\sup_{D_{1/2}}|\nabla h_n|\leq C|x_n|\sup_{D_1} |h_n(y)-m_n|\\
    &\leq C|x_n||h_n(x_n)-m_n|\leq-C|x_n|\log|x_n|\to0, \label{estimategradienthn}
\end{align}
for $n$ sufficiently large. On the other hand, since $T|x_n|\leq \tfrac{s_n}{2}$ and $|x_n|\leq\tfrac{s_n}{2}$, in view of  \eqref{controlhnePn} and \eqref{refinedgradientestimate} we deduce that
\begin{align}
    \nonumber
    |x_n|\sup_{|x|\leq T|x_n|}|\nabla P_n|&\leq C|x_n|\int_{D_2\setminus D_{s_n}}\sup_{|x|\leq T|x_n|}|\nabla_x G(x,y)||y|^{2\alpha}K_n(y)e^{\displaystyle u_n(y)}\,dy\\ \nonumber
    &\leq\frac{C|x_n|}{s_n\log(2/s_n)}\int_{D_2\setminus D_{s_n}}G(x_n,y)|y|^{2\alpha}K_n(y)e^{\displaystyle u_n(y)}\,dy\\ \nonumber
    &=\frac{C|x_n|}{s_n\log(2/s_n)}P_n(x_n)\\
    &\leq \frac{C|x_n|\log(1/|x_n|)}{s_n\log(2/s_n)}=\frac{C\log(1/|x_n|)}{R_n\log(2/s_n)}\leq\frac{C}{R_n}\to0.\label{estimategradientPn}
\end{align}
Therefore, summing the contributions in \eqref{estimategradienthn} and \eqref{estimategradientPn}, we obtain that, for any $T>0$ fixed,
\begin{equation}
    \label{estimategradientHn}
    \sup_{|y|\leq T}|\nabla_y[H_n(|x_n|y)]|
 \leq C\left(|x_n|\log(1/|x_n|)+R_n^{-1}\right)\to0.
\end{equation}

\begin{remark}\label{rem:Hk}
    We notice that the refinement obtained in Lemma \ref{lemmatecnico} is crucial in the above argument since a straightforward application of the standard Harnack inequality to $H_n-m_n$ would just imply that 
    $$|x_n|\sup_{|x|\leq T|x_n|}|\nabla H_n(x)|\leq C\tfrac{-\log|x_n|}{R_n}$$
    and the right-hand side is not guaranteed to converge to zero.
\end{remark}

At this point we define 
\begin{equation}
    \label{rescaledfunction}
v_n(y)=u_n(|x_n|y)+2(1+\alpha)\log|x_n|,\qquad y\in D_{4/|x_n|}(0),
\end{equation}
which satisfies 
\begin{align}
    \label{vnproperties}
    \begin{cases}
        -\lapl v_n=W_ne^{\displaystyle v_n},\qquad\qquad y\in D_{4/|x_n|}(0),\\
W_n(y)=|y|^{2\alpha}K_n(|x_n|y),\\
v_n(\tfrac{x_n}{|x_n|})=2(1+\alpha)\log(\tfrac{|x_n|}{\delta_n})\to+\infty,
    \end{cases}
\end{align}
where
$$0<a\leq K_n(|x_n|y)\leq b<+\infty\quad\text{and}\quad\|\nabla (K_n(|x_n|y))\|_{\infty}\leq |x_n|L\to0.$$
Therefore, possibly along a subsequence, $K_n(|x_n|y)\to c_0$ uniformly on compact sets, with $a\leq c_0\leq b$.
Moreover, from \eqref{massinDsn}, we deduce that
\begin{equation}
    \label{massDRn}
    \int_{D_{R_n}(0)}W_ne^{\displaystyle v_n}<8\pi(1+\alpha).
\end{equation}
Remark that, possibly along a subsequence, $\tfrac{x_n}{|x_n|}\to p_1$, with $|p_1|=1$ and, therefore, $p_1$ is a blow up point for $v_n$. By the concentration-compactness alternative in any fixed ball (\cite{bt}) and using a diagonal argument, we deduce that there exists a nonempty finite set of blow up points $S=\{p_1,\dots,p_m\}$ such that
\[
v_n\to-\infty \qquad \text{uniformly on compact sets of} \,\,\mathbb R^2\setminus S,
\]
\begin{equation}
    \label{concentrationvn}
W_ne^{\displaystyle v_n}\rightharpoonup \sum_{i=1}^m\beta_i\delta_{p_i}, \quad\text{with}\,\, \beta_i\geq8\pi, \qquad\sum_{i=1}^m\beta_i\leq8\pi(1+\alpha).
\end{equation}

\begin{lemma}
\label{lemmaconvergenzariscalate}
Let $y^\ast\notin S$. Then,
\begin{equation*}
    v_n(y)-v_n(y^\ast)\to\Psi(y)-\Psi(y^\ast)\quad\text{ in }\quad\, C^{1}_{loc}(\mathbb{R}^2\setminus S),
\end{equation*}
where
\[
\Psi(y)=-\tfrac{1}{2\pi}\sum_{i=1}^m\beta_i\log|y-p_i|.
\]
In particular, for every compact subset $\omega\subset\mathbb{R}^2\setminus S$, there exists a constant $C_\omega>0$, independent of $n$, such that
\begin{equation}
    \label{boundaryoscillationvn}
    |v_n(x)-v_n(y)|\leq C_\omega,\qquad \text{for every}\,\,x,y\in \omega.
\end{equation}
\end{lemma}
\begin{proof}
By using \eqref{Greenrepresentation2} and rescaling with $|x_n|$, we obtain
\begin{align*}
 v_n(y)&=H_n(|x_n|y)+2(1+\alpha)\log |x_n|\\
&\hspace{2cm}+\frac1{2\pi}\int_{D_{R_n}}
 \log\left(\frac{|4-|x_n|^2\overline\eta y|}{2|x_n||y-\eta|}\right)W_ne^{\displaystyle v_n}\,d\eta\\
 &=A_n(y)+B_n(y)
\end{align*}
where
\[
A_n(y):=H_n(|x_n|y)+2(1+\alpha)\log |x_n|+\frac1{2\pi}\int_{D_{R_n}}
 \log\left(\frac{|4-|x_n|^2\overline\eta y|}{2|x_n|}\right)W_ne^{\displaystyle v_n}\,d\eta,
\]
\[
B_n(y):=-\frac{1}{2\pi}\int_{D_{R_n}}
 \log{|y-\eta|}W_ne^{\displaystyle v_n}\,d\eta,
\]
Let $T>0$. For $|y|\leq T$ and $|\eta|\leq R_n$, we observe that $|x_n|^2|\eta||y|\leq T|x_n|^2R_n=T|x_n|s_n\to0$ and, therefore,
\[
\sup_{|y|\leq T}\left|\nabla_y\int_{D_{R_n}}
 \log\left(\frac{|4-|x_n|^2\overline\eta y|}{2|x_n|}\right)W_ne^{\displaystyle v_n}\,d\eta\right|\leq C|x_n|^2R_n\int_{D_{R_n}}W_ne^{\displaystyle v_n}\,d\eta\leq C|x_n|s_n\to0.
\]
Combining this result with \eqref{estimategradientHn}, we obtain that 
\[
\|\nabla A_n\|_{L^\infty(D_T)}\to0
\]
and in particular, considering $T>|y^\ast|$, by the mean value theorem, we have that
\begin{equation}
    \label{convergenceregularpart}
    A_n(y)-A_n(y^\ast)\to0,\quad \text{ in }\quad C^1(\bar D_T).
\end{equation}
We are left to show that 
\begin{equation}
\label{convergence B_n}
Q_n:=B_n-B_n(y^\ast)\to\Psi-\Psi(y^\ast) \quad\text{ in }\quad C^1_{\mathrm{loc}}(\mathbb R^2\setminus S)
\end{equation}
Firstly, we observe that
\[
 Q_n(y)=-\frac1{2\pi}\int_{D_{R_n}}
 \log\left(\frac{|y-\eta|}{|y^\ast-\eta|}\right)W_ne^{\displaystyle v_n}\,d\eta
\]
and
\[
\nabla Q_n(y)=-\frac1{2\pi}\int_{D_{R_n}}
 \frac{y-\eta}{|y-\eta|^2}W_ne^{\displaystyle v_n}\,d\eta.
\]
Let $\omega\subset\mathbb R^2\setminus S$ be a compact set that contains $y^\ast$ and choose $T>0$ such that $\omega\cup S\subset D_T(0)$. Set $A>2T$ and $\rho$ sufficiently small such that $\rho<\tfrac14\operatorname{dist}(\omega,S)$ and $B_{\rho} (p_i)$ are all disjoint. We decompose the whole domain into
\[
D_{R_n}(0)= U_\rho\cup E_{\rho,A}\cup F_{n,A},
\]
with
\[
U_\rho:=\bigcup_{i=1}^m B_{\rho}(p_i),\quad E_{\rho,A}:=D_A(0)\setminus U_\rho,\quad F_{n,A}=D_{R_n}(0)\setminus D_A(0).
\]
Since $\omega$ is at a positive distance from $S$, for every $y\in \omega$ and $\eta\in B_\rho(p_i)$ we have
\[
    |y-\eta|\geq \frac{3}{4}\operatorname{dist}(\omega,S).
\]
Consequently, the kernels inside $Q_n$ and $\nabla Q_n$ are uniformly
Lipschitz continuous with respect to $\eta$ in these regions. In particular,
\[
    \left|\log\left(\frac{|y-\eta|}{|y^\ast-\eta|}\right)-\log\left(\frac{|y-p_i|}{|y^\ast-p_i|}\right)\right|=\left|\log\left(\frac{|y-\eta||y^\ast-p_i|}{|y^\ast-\eta||y-p_i|}\right)\right|\leq C_\omega\rho
\]
and
\[
    \left|\frac{y-\eta}{|y-\eta|^2}-\frac{y-p_i}{|y-p_i|^2}\right|\leq C_\omega\rho,
\]
for some positive constant $C_\omega$ which depends only on $\omega$.
It follows that
\begin{align*}
&\left|\int_{B_\rho(p_i)}\log\left(\frac{|y-\eta|}{|y^\ast-\eta|}\right)W_ne^{\displaystyle v_n}\,d\eta
    -\beta_i \log\left(\frac{|y-p_i|}{|y^\ast-p_i|}\right)
\right|\\
&\hspace{2cm}\qquad\leq
C_\omega\rho\,\int_{B_\rho(p_i)}W_ne^{\displaystyle v_n}\,d\eta+\left|\log\left(\frac{|y-p_i|}{|y^\ast-p_i|}\right)\right|\left|\int_{B_\rho(p_i)}W_ne^{\displaystyle v_n}\,d\eta-\beta_i\right|.
\end{align*}
Summing over $i$ and using the uniform mass bound and \eqref{concentrationvn},
we obtain
\begin{equation}
\label{eq:bubbles-potential}
\limsup_{n\to\infty}\sup_{y\in \omega}
\left|
    \int_{U_\rho}\log\left(\frac{|y-\eta|}{|y^\ast-\eta|}\right)W_ne^{\displaystyle v_n}\,d\eta
    -\sum_{i=1}^{m}\beta_i \log\left(\frac{|y-p_i|}{|y^\ast-p_i|}\right)
\right|\leq C_\omega\rho.
\end{equation}
The same argument works for the derivative and we obtain
\begin{equation}
\label{eq:bolle-gradiente}
\limsup_{n\to\infty}\sup_{y\in \omega}
\left|
\int_{U_\rho}\frac{y-\eta}{|y-\eta|^2}W_ne^{\displaystyle v_n}\,d\eta
    -\sum_{i=1}^{m}\beta_i \frac{y-p_i}{|y-p_i|^2}
\right|\leq C_\omega\rho.
\end{equation}
Let us consider now $\eta\in E_{\rho,A}$, then, since $\bar E_{\rho,A}\subset\mathbb R^2\setminus S$ and by \eqref{concentrationvn}, we have  $\|W_ne^{\displaystyle v_n}\|_{L^\infty(E_{\rho,A})}\to 0$ and, consequently, we obtain
\begin{equation}
\label{eq:regular-region-potential}
\begin{aligned}
\sup_{y\in \omega}
\left|\int_{E_{\rho,A}}
    \log\left(\frac{|y-\eta|}{|y^\ast-\eta|}\right)W_ne^{\displaystyle v_n}\,d\eta\right|&\leq \|W_ne^{\displaystyle v_n}\|_{L^\infty(E_{\rho,A})}\sup_{y\in \omega}\int_{E_{\rho,A}}|\log|y-\eta|-\log|y^\ast-\eta||\,d\eta
\\
&\leq
    2\|W_ne^{\displaystyle v_n}\|_{L^\infty(E_{\rho,A})}
    \int_{D_{A+T}(0)}|\log|z||\,dz\\
&=
    4\pi\|W_ne^{\displaystyle v_n}\|_{L^\infty(E_{\rho,A})}
    \int_0^{A+T}r|\log r|\,dr
    \longrightarrow 0
\end{aligned}
\end{equation}
and
\begin{equation}
\label{eq:regular-region-gradient}
\begin{aligned}
\sup_{y\in \omega}
\left|\int_{E_{\rho,A}}\frac{y-\eta}{|y-\eta|^2}W_ne^{\displaystyle v_n}\,d\eta\right|&\leq
    \|W_ne^{\displaystyle v_n}\|_{L^\infty(E_{\rho,A})}\sup_{y\in \omega}
    \int_{D_A}\frac{d\eta}{|y-\eta|}\\
&\leq\|W_ne^{\displaystyle v_n}\|_{L^\infty(E_{\rho,A})}
    \int_{D_{A+T}(0)}\frac{dz}{|z|}
\\
&=2\pi(A+T)\|W_ne^{\displaystyle v_n}\|_{L^\infty(E_{\rho,A})}
    \longrightarrow 0.
\end{aligned}
\end{equation}
Finally, let $\eta\in F_{n,A}$. Concerning the normalized logarithmic potential, we
observe that
\[
\begin{aligned}
\log\left(\frac{|y-\eta|}{|y^\ast-\eta|}\right)=\int_0^1
\frac{
    (y^\ast+t(y-y^\ast)-\eta)\cdot(y-y^\ast)
}{
    |y^\ast+t(y-y^\ast)-\eta|^2
}\,dt.
\end{aligned}
\]
Since the line segment joining $y^\ast$ and $y$ is
contained in $D_T(0)$, we have
\[
    |y^\ast+t(y-y^\ast)-\eta|
    \geq |\eta|-T
    \geq \frac{|\eta|}{2},
    \qquad t\in[0,1].
\]
Therefore,
\[
    \left|\log\left(\frac{|y-\eta|}{|y^\ast-\eta|}\right)\right|
    \leq\frac{2|y-y^\ast|}{|\eta|}
    \leq\frac{4T}{A}.
\]
It follows that
\begin{equation}
\label{eq:tail-potential}
\begin{aligned}
&\sup_{y\in \omega}
\left|
    \int_{F_{n,A}}
    \log\left(\frac{|y-\eta|}{|y^\ast-\eta|}\right)W_ne^{\displaystyle v_n}\,d\eta\right|\leq
    \frac{4T}{A}\int_{D_{R_n}(0)}W_ne^{\displaystyle v_n}\,d\eta
    \leq\frac{32\pi(1+\alpha) T}{A}.
\end{aligned}
\end{equation}
Concerning the derivative, for $|\eta|>A>2T$ and $y\in \omega$, we have
\[
    |y-\eta|\geq |\eta|-T
    \geq\frac{|\eta|}{2}.
\]
Hence,
\[
    \left|\frac{y-\eta}{|y-\eta|^2}\right|
    \leq\frac{2}{|\eta|}
    \leq\frac{2}{A}.
\]
Using once again the uniform mass bound, we obtain
\begin{equation}
\label{eq:tail-gradient}
\begin{aligned}
&\sup_{y\in \omega}
\left|
    \int_{F_{n,A}}
    \frac{y-\eta}{|y-\eta|^2}W_ne^{\displaystyle v_n}\,d\eta\right|
\leq
    \frac{2}{A}\int_{D_{R_n}(0)}W_ne^{\displaystyle v_n}\,d\eta
    \leq\frac{16\pi(1+\alpha)}{A}.
\end{aligned}
\end{equation}

Combining \eqref{eq:bubbles-potential},
\eqref{eq:regular-region-potential}
and \eqref{eq:tail-potential}, we obtain
\[
\limsup_{n\to\infty}\sup_{y\in \omega}
\left|
    \int_{D_{R_n}(0)} \log\left(\frac{|y-\eta|}{|y^\ast-\eta|}\right)W_ne^{\displaystyle v_n}\,d\eta-\sum_{i=1}^{m}\beta_i  \log\left(\frac{|y-p_i|}{|y^\ast-p_i|}\right)
\right|\leq C_\omega\rho+\frac{C_\omega}{A}.
\]
Similarly, by \eqref{eq:bolle-gradiente},
\eqref{eq:regular-region-gradient}
and \eqref{eq:tail-gradient}, we have
\[
\limsup_{n\to\infty}\sup_{y\in \omega}
\left|
    \int_{D_{R_n}(0)}\frac{y-\eta}{|y-\eta|^2}W_ne^{\displaystyle v_n}\,d\eta-\sum_{i=1}^{m}\beta_i  \frac{y-p_i}{|y-p_i|^2}
\right|\leq C_\omega\rho+\frac{C_\omega}{A}.
\]
Letting first $\rho\to0$ and subsequently
$A\to+\infty$, we conclude that
\[
    Q_n(y)\longrightarrow
    -\frac{1}{2\pi}
    \sum_{i=1}^{m}\beta_i
    \log\frac{|y-p_i|}{|y^\ast-p_i|}\,\,\,\,\text{and}\,\,\,\, \nabla Q_n(y)\longrightarrow
    -\frac{1}{2\pi}
    \sum_{i=1}^{m}\beta_i
    \frac{y-p_i}{|y-p_i|^2},
\]
uniformly on $\omega$. Since $\omega$ was arbitrary, we prove our result. Finally, for any compact subset
$\omega\subset\mathbb R^2\setminus S$, let us set
\[
r_n(y):=v_n(y)-v_n(y^\ast)
-\Psi(y)+\Psi(y^\ast).
\]
Then, for every $x,y\in \omega$, we have that
\[
|v_n(x)-v_n(y)|
\leq |\Psi(x)-\Psi(y)|
+|r_n(x)|+|r_n(y)|
\leq \operatorname{osc}_\omega\Psi
+2\|r_n\|_{L^\infty(\omega)}
\leq C_\omega,
\]
which proves \eqref{boundaryoscillationvn}.
\end{proof}

We observe that Lemma \ref{lemmaconvergenzariscalate} provides a uniform bound on the oscillation of the rescaled solutions $v_n$ on compact subsets away from their blow up set, thereby recovering the property established in Lemma 2.2 of \cite{bcct}, which was obtained therein assuming a uniform bound on the boundary oscillation of the original solutions $u_n$ and of the global mass. At this point, the remaining part of the proof follows step by step the moving plane argument worked out in the proof of Theorem 1.3 of \cite{bcct}. We provide the details in the Appendix for the sake of completeness.
In particular, we obtain
\begin{equation}
    \label{contradiction2}
    u_n(x_n)+\inf_{D_4}u_n\leq C
\end{equation}
for some constant $C$ independent of $n$, thus contradicting \eqref{contradiction}.
\end{proof}

\smallskip

\section{appendix}
We follow the argument in \cite{bcct}. Firstly, we show that the blow up set of $v_n$ does not contain $0$. Suppose, by contradiction, that $0\in S$. Then, since $v_n$ satisfies \eqref{boundaryoscillationvn}, we apply the quantization result in \cite{bt} to deduce that
\[
\int_{D_\delta(0)}
W_ne^{\displaystyle v_n}\,dx
=8\pi(1+\alpha)+o(1),
\qquad \text{as } n\to+\infty,
\]
for $\delta>0$ sufficiently small.
On the other hand, $p_1$ is a regular blow up point for $v_n$, therefore, applying \cite{yy}, we deduce that
\[
\int_{D_\delta(p_1)}
W_ne^{\displaystyle v_n}\,dx
=8\pi+o(1),
\qquad \text{as } n\to+\infty,
\]
for $\delta>0$ sufficiently small.
Combining these facts, together with \eqref{massDRn}, we get 
\[
8\pi(1+\alpha)\geq\int_{D_{R_n}(0)}W_ne^{\displaystyle v_n}\geq \int_{D_\delta(p_1)}
W_ne^{\displaystyle v_n}+\int_{D_\delta(0)}
W_ne^{\displaystyle v_n}
=8\pi(1+\alpha)+8\pi+o(1),
\]
which is a contradiction. \\
Therefore, the blow up set of $v_n$ is $S=\{p_1,\ldots,p_m\}
\subset\mathbb{R}^2\setminus\{0\}$.
We are now in position to apply Li's result (\cite{yy}) and conclude 
\begin{equation}
\label{Liprofile}
f_{n,i}(y):=
\left|
v_n(y)-
\log\left(\frac{e^{v_n(y_{n,i})}}
{\left(
1+\sigma_{n,i}e^{v_n(y_{n,i})}
|y-y_{n,i}|^2
\right)^2}\right)
\right|
\leq C,
\end{equation}
for any $y\in\overline{D}_{r_0}(p_i)$,
where $y_{n,i}\in\overline{D}_{r_0}(p_i)$ satisfies
\[
v_n(y_{n,i})
=\max_{\bar D_{r_0}(p_i)}v_n,
\qquad \text{and}\qquad y_{n,i}\to p_i,
\]
\[
\sigma_{n,i}
=\frac{1}{8}|y_{n,i}|^{2\alpha}
K_n(|x_n|y_{n,i})
\to
\frac{1}{8}|p_i|^{2\alpha}c_0
:=\sigma_i>0,
\qquad \text{as } n\to+\infty,
\]
for $i\in\{1,\ldots,m\}$ and $r_0>0$
sufficiently small. In particular, since $x_n$ is a maximum point
of $u_n$ in $\overline D_1$, we may choose
$y_{n,1}=x_n/|x_n|$.
In addition, for any fixed $R>0$,
$f_{n,i}$ defined in \eqref{Liprofile} satisfies
\begin{equation}
\label{profilelocal}
\max\left\{
f_{n,i}(y)\,|\,
e^{v_n(y_{n,i})}|y-y_{n,i}|^2\leq R
\right\}
\to 0,
\qquad \text{as } n\to+\infty,
\end{equation}
for $i\in\{1,\ldots,m\}$. Note in particular that there exists $\tau_0>0$
such that
\begin{equation}
\label{boundfrombelowmaxima}
|y_{n,i}|\geq\tau_0,
\qquad
\forall i\in\{1,\ldots,m\},
\end{equation}
for $n$ sufficiently large.
Set $\mathcal D_i=D_{r_0}(p_i)$ and take $R_0>0$
sufficiently large so that
$\bigcup_{i=1}^{m}\mathcal D_i\Subset D_{R_0}$
and let
$\Omega_1=D_{R_0}\setminus
\overline{\bigcup_{i=1}^{m}\mathcal D_i}$.
By \eqref{concentrationvn}, $v_n\to-\infty$ on $\overline{\Omega}_1$ and, by \eqref{boundaryoscillationvn}, we conclude
\begin{equation}
\label{oscillationcontrolled}
\max_{\overline{\Omega}_1}v_n
-\min_{\overline{\Omega}_1}v_n
\leq C_0,
\end{equation}
for a suitable constant $C_0>0$.
Hence, by combining \eqref{Liprofile} and \eqref{oscillationcontrolled} we have
\[
\left|
v_n\left(\frac{x_n}{|x_n|}\right)
-v_n(y_{n,i})
\right|
=O(1),
\qquad \forall i\in\{1,\ldots,m\}.
\]

Thus, let
\begin{equation}
\label{rescaling}
\mu_n
=e^{v_n\left(\frac{x_n}{|x_n|}\right)}
=|x_n|^{2(1+\alpha)}e^{u_n(x_n)}\to+\infty,
\end{equation}
then we derive
\begin{equation}
\label{finalprofile}
\left|
v_n(y)-
\log
\left(\frac{\mu_n}
{\left(
1+\sigma_0\mu_n
\min_{i\in\{1,\ldots,m\}}|y-y_{n,i}|^2
\right)^2}\right)
\right|
\leq C,
\quad \forall y\in D_{R_0},
\end{equation}
for suitable constants $\sigma_0>0$ and $C>0$.

Let
\[
Q=\{(t,\theta):t\leq\log 2,\ 
\theta\in[0,2\pi]\},
\]
\[
w_n(t,\theta)
=u_n(e^t\cos\theta,e^t\sin\theta)
+2(1+\alpha)t-\frac{L}{a}e^t,
\qquad \text{in } Q.
\]

Hence $w_n$ satisfies
\[
-\Delta w_n
=\widetilde{V}_n(t,\theta)e^{w_n}
+\frac{L}{a}e^t
\qquad \text{in } Q,
\]
with
\[
\widetilde{V}_n(t,\theta)
=e^{\frac{L}{a}e^t}
K_{n}(e^t\cos\theta,e^t\sin\theta).
\]
As in \cite{bls}, we deduce that
\begin{equation}
\label{monotonicity}
\frac{\partial}{\partial t}
\left\{
\widetilde{V}_n(t,\theta)
e^\xi
+\frac{L}{a}e^t
\right\}
\geq 0,
\qquad \text{in } Q,
\end{equation}
for every $\xi\in\mathbb{R}$ and that there exists
$\lambda_0<0$ depending on $n$, such that
\begin{equation*}
w_n(2\lambda-t,\theta)-w_n(t,\theta)<0,
\qquad
\forall\lambda\leq\lambda_0,\quad\forall \lambda<t\leq\log 2
\quad \text{and}\quad
\theta\in[0,2\pi].
\end{equation*}

Hence, the following quantity is well defined
\[
\lambda_n
=\sup\{\lambda<\log 2 \,|\, w_n(2\mu-t,\theta)-w_n(t,\theta)<0,
\,\,\forall\mu\leq\lambda,\,\,
\forall\mu<t\leq\log 2
\,\, \text{and}\,\,
\theta\in[0,2\pi] \}.
\]

\medskip

Now, we prove that there exists $C>0$ such that
\begin{equation}
\label{boundlambdan}
\lambda_n\leq\log|x_n|+C.
\end{equation}

To establish \eqref{boundlambdan}, let
$\theta_n\in[0,2\pi)$ such that
$e^{i\theta_n}=\frac{x_n}{|x_n|}$. By means of \eqref{profilelocal}, we notice that,
\begin{equation}
\label{comparingprofile}
w_n(\log|x_n|+s,\theta_n)=
\log \mu_n+2(1+\alpha)s-2\log\left(
1+\sigma_1\mu_n(1-e^s)^2
\right)+o(1),
\end{equation}
uniformly in $|s|\leq\tfrac{1}{\sqrt{\mu_n}}$,
with $\sigma_1=\tfrac{c_0}{8}$. Let $\epsilon_n=\log\left(1+\frac{1}{\sqrt{\mu_n}}\right)\to0$, then $\mu_n(1-e^{\epsilon_n})^2=1$ and, from \eqref{comparingprofile} we deduce that
\begin{align*}
& w_n(\log|x_n|,\theta_n)
-w_n(\log|x_n|+\epsilon_n,\theta_n)=2\log(1+\sigma_1)-2(1+\alpha)\epsilon_n+o(1)\to 2\log(1+\sigma_1)>0,
\end{align*}
as $n\to+\infty$. Therefore, by choosing
\[
\lambda
=\log|x_n|
+\frac{1}{2}\log\left(
1+\frac{1}{\sqrt{\mu_n}}
\right)
\qquad \text{and}\qquad
t=\log|x_n|
+\log\left(
1+\frac{1}{\sqrt{\mu_n}}
\right),
\]
we have that
\[
\begin{aligned}
w_n(2\lambda-t,\theta_n)-w_n(t,\theta_n)
&=w_n(\log|x_n|,\theta_n)
-w_n(\log|x_n|+\epsilon_n,\theta_n)>0,
\end{aligned}
\]
for $n$ large and \eqref{boundlambdan} is established.

Moreover, by using the maximality of $\lambda_n$ together with \eqref{monotonicity} and the maximum principle, we deduce that,
\begin{equation}
\label{boundwn}
\min_{\theta\in[0,2\pi)}
w_n(\log 2,\theta)
\leq
\max_{\theta\in[0,2\pi)}
w_n(2\lambda_n-\log2,\theta).
\end{equation}
From the definition of $w_n$, we have that
\[
w_n(2\lambda_n-\log2,\theta)=v_n(\tfrac{e^{2\lambda_n}e^{i\theta}}{2|x_n|})+2(1+\alpha)\log(\tfrac{e^{2\lambda_n}}{2|x_n|})-\tfrac{L}{2a}e^{2\lambda_n}.
\]
By using \eqref{boundlambdan}, $\tfrac{e^{2\lambda_n}e^{i\theta}}{2|x_n|}\to0$, uniformly in $\theta$, while by \eqref{boundfrombelowmaxima} and \eqref{finalprofile} we deduce that
\[
v_n(\tfrac{e^{2\lambda_n}e^{i\theta}}{2|x_n|})\leq-\log\mu_n+C.
\]
Therefore,
\begin{equation}
\label{finalestimate}
\begin{aligned}
w_n(2\lambda_n-\log2,\theta)
&\leq v_n(\tfrac{e^{2\lambda_n}e^{i\theta}}{2|x_n|})+2(1+\alpha)(2\lambda_n-\log2-\log|x_n|)\\
&\leq-\log\mu_n+2(1+\alpha)\log|x_n|+C\\
&=-u_n(x_n)+C.
\end{aligned}
\end{equation}
Since
\[
\inf_{D_4}u_n+2(1+\alpha)\log2-2\tfrac La\leq \min_{\theta\in[0,2\pi)}w_n(\log2,\theta),
\]
combining with \eqref{boundwn} and \eqref{finalestimate}, we get \eqref{contradiction2}.

\section*{Data availability statement}
Data sharing is not applicable to this article as no datasets were generated or analysed during the current study.

\section*{Conflict-of-Interest Statement}
The authors declare that they have no conflict of interest.

\section*{Declaration of AI use}
The authors have been investigating the problems addressed in this paper for several years. During the development of this work, the authors used AI tools to test new approaches based on ideas in \cite{B4} and \cite{bcct}, discuss selected technical aspects, conduct literature searches and improve the exposition. In particular, AI suggested the relevance of the argument in reference [33], which was previously unknown to the authors and pointed out that the oscillation of the function $H_n$ in (5.7) is small at the scale $|x_n|$; see Remark 5.2. All other proof strategies and mathematical arguments were developed by the authors, who take full responsibility for the entire content of the paper.


\begin{thebibliography}{99}

%\bibitem{ALN} A. Ataei, D. Lundholm, D.-T. Nguyen \emph{A generalized liouville equation and magnetic stability}, arXiv:2404.09332v2

\bibitem{B1} D. Bartolucci, {\em A ``Sup + C Inf" inequality for the
equation $-\Delta u=\frac{V}{|x|^{2\alpha}} e^{u}$},
{Proc. Royal Soc. of Edinburgh}, {\bf 140A} (2010), 1119-1139.

\bibitem{B2} D. Bartolucci, {\em A ``Sup + Inf" inequality for Liouville
type equations with weights }, Jour.  d'Analyse
Mathematique, {\bf 117} (2012), 29-46; DOI:10.1007/s11854-012-0013-7.

\bibitem{B3} D. Bartolucci, {\em A ``Sup + C Inf" inequality for Liouville
type equations with singular potentials },
Mathematische Nachrichten, {\bf 284} (2011), 1639-1651.

\bibitem{B4} D. Bartolucci, \emph{A sup$\times$inf-type inequality for conformal metrics on Riemann surfaces with conical singularities},
Journal of Mathematical Analysis and Applications,
{\bf 403}  (2013), 571-579.


\bibitem{bcct} D. Bartolucci, C.C. Chen, C.S. Lin and G. Tarantello,
{\em Profile of Blow Up Solutions To Mean Field Equations with Singular Data},
Comm. in P. D. E.,  {\bf 29}(7-8) (2004), 1241--1265.

\bibitem{BCW0} 
D. Bartolucci, P. Cosentino and L. Wu,
{\em A Harnack-Type Inequality for a Perturbed Singular Liouville Equation},
Math. Nachr. {\bf 299} (2026), no. 6, 1398-1417, doi: 10.1002/mana.70149.

\bibitem{bart-5} D. Bartolucci, A. Jevnikar, Y. Lee, W. Yang, \emph{Non degeneracy, Mean Field Equations and the Onsager theory of 2D turbulence}, {Arch. Ration. Mech. An.},  \textbf{230} (2018), 397-426.

\bibitem{bt}
D. Bartolucci and G. Tarantello, {\em Liouville type equations with singular data and their applications to periodic multivortices for the electroweak theory}, Commun. Math. Phys. {\bf 229}, 3-47 (2002).

\bibitem{bt2} D. Bartolucci and G. Tarantello, {\em The Liouville Equation with Singular Data:
A Concentration-Compactness Principle via a Local
Representation Formula}, J. Diff. Equations {\bf 185} (2002), 161--180.

\bibitem{bart-6} D. Bartolucci, W. Yang, L. Zhang, 
\emph{Asymptotic analysis and uniqueness of blowup solutions of
non-quantized singular mean field equations},
Math. Ann., \textbf{395} (2026), 103.

\bibitem{bls}
H. Brezis, Y.Y. Li and I.Shafrir, {\em A sup+inf inequality for some nonlinear elliptic equations involving exponential nonlinearities}, Jour. Func. Analysis, {\bf 115} (1993), 344--358.

\bibitem{bm}
H. Brezis and F. Merle,
{\em Uniform estimates and blow-up behaviour for
solutions of $-\Delta u = V(x)e^{u}$ in two dimensions},
{Comm. in P.D.E.,}  {\bf 16}(8,9) (1991), 1223--1253.

\bibitem{cl} C. C. Chen and C. S. Lin, {\em A sharp sup+inf inequality for a nonlinear elliptic
equation in $\rdue$}, Comm. An. Geom.,  {\bf 6}(1) (1998), 1--19.

\bibitem{cl-1} C. C. Chen and C. S. Lin, {\em Estimates of the conformal scalar curvature equation via the
method of moving planes}, Comm. Pure Appl. Math.,  {\bf 50} (1997), 971-1017.

%\bibitem{cl-2} C.C. Chen, C.S. Lin, %Topological degree for a mean field %equation on Riemann surfaces. \emph{Comm. %Pure Appl. Math.}, \textbf{56} (2003), %1667-1727.
		
\bibitem{cl-3} C.C. Chen, C.S. Lin, \emph{Mean field equation of Liouville type with singular data: topological degree}, {Comm. Pure Appl. Math.}, \textbf{68} (2015), 6, 887-947.

\bibitem{CW}
K.S. Chou and T. Y. H. Wan,
\emph{Asymptotic radial symmetry for solutions of $\Delta u + \exp{u}=0$
in a punctured disc}, {Pacific J. of Math.} {\bf 163}(2) (1994), 269-276.

\bibitem{Cos} P. Cosentino, \emph{A Harnack type inequality for singular Liouville type equations},
{Journal of Functional Analysis}, {\bf 289}, (2025) 111003.

\bibitem{dwz}
T. D'Aprile, J.C. Wei, L. Zhang, {\em On non-simple blowup solutions of Liouville equation}, American. J.  Math., {\bf 148} (2026), 1329-1369.

\bibitem{FKR} R. Fournier, D. Kraus, O. Roth, \emph{A Schwarz lemma for locally univalent meromorphic functions},
Proc.  A. M. S. {\bf 148} (2020), 3859-3870.

\bibitem{GKR} J. Grahl, D. Kraus, O. Roth, \emph{Blow-up solutions of Liouville's equation and quasi-normality},
Computational Methods and Function Theory {\bf 20} (2020), 677-692.

\bibitem{KN} S. Kobayashi, K. Nomizu, {\em Foundations of Differential Geometry}, vol. 1,
Interscience Publishers, John Wiley and Sons, 1963.

\bibitem{llty}
Y. Lee, C.S. Lin, G. Tarantello and W. Yang, {\em Sharp estimates for solutions of mean
field equations with collapsing singularity}. Comm. Partial Differ. Eq. {\bf 42} (2017), 1549-1597.

\bibitem{yy} Y.Y. Li,  {\em Harnack type inequality: the method of moving planes},
Comm. Math. Phys.  {\bf 200} (1999), 421-444.

\bibitem{Li-1} Y.Y. Li, \emph{Prescribing scalar curvature on $\mathbb{S}^n$ and related problems, Part I},
J. Diff. Equations {\bf 120} (1995), 319-410.

\bibitem{Li-2} Y.Y. Li, \emph{Prescribing scalar curvature on $\mathbb{S}^n$ and related problems, Part II: Existence and
compactness}, Comm. Pure Appl. Math. {\bf 49} (1996), 541-597.

\bibitem{LiZ0}
 Y.Y. Li, L. Zhang, \emph{Liouville-type theorems and harnack-type inequalities for semilinear elliptic equations}, J. Anal. Math. {\bf 90}, 27-87 (2003).

\bibitem{LiZ} Y.Y. Li, L. Zhang, \emph{A Harnack type inequality for the Yamabe equation in low dimensions}, Calc. Var. P.D.E. {\bf 20} (2004), 133-151.

\bibitem{Lin-0} C.S. Lin, \emph{Spherical metrics with one singularity and odd integer angle on flat tori, I}, Jour. Diff. Geom. {\bf 132} (2026), 
321-382.

\bibitem{Lio}
J. Liouville,
{\em Sur L' \'Equation aux diff\'erences partielles
$\frac{d^{2} \log{\lm}}{du dv} \pm \frac{\lm}{2 a^{2}}=0$},\\
{Jour. Math. Pure App.}, {\bf 18} (1853), 71-72.

\bibitem{Pic} E. Picard,
{\em De l'int\'{e}gration de l'\'{e}quation $\Delta u= e^u$ sur une surface de Riemann ferm\'{e}e},
J. Crelle, {\bf 130}, (1905), 243-258.

\bibitem{Shoen} Schoen, R.: Courses at Stanford University (1988) and New York University (1989)
(unpublished).

\bibitem{S} I.Shafrir, \emph{A Sup+Inf inequality for the equation $-\Delta u = V(x)e^{u}$},
{C. R. Acad. Sci. Paris},  {\bf 315} (1992), 159--164.


\bibitem{St} N. Steinmetz, \emph{Normal families and linear differential equations},
Journal d'Analyse Math\'ematique, {\bf 117} (2012), 129-132.

\bibitem{T1} G. Tarantello,
\emph{ A Harnack inequality for Liouville type equations with singular sources},
Indiana. Univ. Math. Jour., {\bf 54}(2) (2005), 599--616.

\bibitem{T3} G. Tarantello,
{\sl Analytical aspects of Liouville type equations with singular sources}, Handbook Diff. Eqs., North Holland,
Amsterdam, Stationary partial differential equations, {\bf I} (2004), 491--592.

\bibitem{T2} G. Tarantello,
{\sl Self-Dual Gauge Field Vortices: an analytical approach}, Progress in Nonlin. Diff. Eqs. and Appl.,
{\bf 72} (2008).

\bibitem{T4} G. Tarantello, \emph{Blow-up analysis for a cosmic strings equation}, \emph{Jour. Funct. An.}, \textbf{272} (2017), 255-338.

\bibitem{T5} G. Tarantello,
 \emph{Asymptotics for minimizers of the Donaldson
functional in Teichmuller theory}, Adv. Math. \textbf{425} (2023), 109090.

\bibitem{wz1}
J.C. Wei, L. Zhang, {\em Estimates for Liouville equation with quantized singularities}, Adv. Math., {\bf 380} (2021), 107606.

%\bibitem{wz2}
%J.C. Wei and L. Zhang, {\em Vanishing estimates %for Liouville equation with quantized %singularities}, Proc. Lond. Math. Soc., {\bf 3} %(2022), 1-26.

%\bibitem{Y} Y. Yang, {\sl Solitons in Field Theory and Nonlinear Analysis},
%Springer Monographs in Mathematics, Springer-Verlag, New-York, 2001

\bibitem{zhang1} L. Zhang, \emph{Blowup solutions of some nonlinear elliptic equations involving
exponential nonlinearities}, {Comm. Math. Phys}, \textbf{268} (2006), 105-133.
		
\bibitem{zhang2} L. Zhang, \emph{Asymptotic behavior of blowup solutions for elliptic equations with exponential nonlinearity and singular data}, {Comm. Contemp. Math.}, \textbf{11} (2009), 395-411.


\end{thebibliography}
\end{document}